\documentclass[11pt]{article}

\usepackage{amssymb}

\begin{document}

\def\pni{\par \noindent}
\def\vsh{\smallskip}
\def\s{\smallskip}
\def\vs{\medskip}
\def\vvs{\bigskip}
\def\vvvs{\bigskip\medskip} 
\def\vsp{\vsh\pni}
\def\vsn{\vsh\pni}
\def\cen{\centerline}
\def\ra{\item{a)\ }} \def\rb{\item{b)\ }}   \def\rc{\item{c)\ }}
\def\eg{{\it e.g.}\ } \def\ie{{\it i.e.}\ }


\def\sg{\hbox{sign}\,}
\def\sgn{\hbox{sign}\,}
\def\sign{\hbox{sign}\,}
\def\e{{\rm e}}
\def\exp{{\rm exp}}
\def\ds{\displaystyle}
\def\dis{\displaystyle}
 \def\q{\quad}    \def\qq{\qquad} 
\def\lan{\langle}\def\ran{\rangle}
\def\l{\left} \def\r{\right}
\def\lt{\left} \def\rt{\right}  
\def\lra{\Longleftrightarrow}
\def\d{\partial}
\def\dr{\partial r}  \def\dt{\partial t}
\def\dx{\partial x}   \def\dy{\partial y}  \def\dz{\partial z}
\def\rec#1{{1\over{#1}}}
\def\zr{z^{-1}}



\def\hatt{\widehat}
\def\epsilons{{\widetilde \epsilon(s)}}
\def\sigmas{{\widetilde \sigma (s)}}
\def\fs{{\widetilde f(s)}}
\def\Js{{\widetilde J(s)}}
\def\Gs{{\widetilde G(s)}}
\def\Fs{{\wiidetilde F(s)}}
 \def\Ls{{\widetilde L(s)}}
\def\L{{\mathcal L}} 
\def\F{{\mathcal F}} 
\def\NN{{\bf N}}
\def\RR{{\bf R}}
\def\CC{{\bf C}}
\def\ZZ{{\bf Z}} 

\begin{center}
{\bf The Mittag-Leffler function in the thinning theory \\ for  renewal processes}

\vs
{Rudolf GORENFLO}$^{(1)}$ and 
 {Francesco MAINARDI}$^{(2)}$
\end{center}

\vsn
$\null^{(1)}$
 First Mathematical Institute,
 Free   University of Berlin, \\
 Arnimallee  3, D-14195 Berlin, Germany. \\
Deceased (1930-2017)
\\ [0.25 truecm]
$\null^{(2)}$
Department of Physics and Astronomy, University of Bologna, and INFN, \\
Via Irnerio 46, I-40126 Bologna, Italy. \\
E-mail: francesco.mainardi@bo.infn.it 
\vvs \vsn

\begin{abstract}
The main purpose of this note is to  point out the relevance of the Mittag-Leffler probability distribution  in the so-called thinning theory for a renewal process with a queue of power law type.
This theory, formerly considered by Gnedenko and Kovalenko in 1968 without the explicit reference to the Mittag-Leffler function, was used by the authors in the theory of continuous random walk and consequently of fractional diffusion in  a plenary lecture by the late Prof Gorenflo at a Seminar on Anomalous Transport held in Bad-Honnef in July 2006, published in a 2008 book.
After recalling the basic theory of renewal processes including the standard and the fractional Poisson processes, here we have revised the original approach by Gnedenko and Kovalenko for convenience of the experts of stochastic processes who are not aware of the relevance of the Mittag-Leffler functions.

\end{abstract}

\noindent
This note is devoted to the memory of the late Professor Rudolf Gorenflo, passed away on 20 October 2017 at the age of 87.
It has been published in the journal {\bf Theory of Probability and Mathematical Statistics, Vol. 98,  No  1, pp. 100-108 (2018)}.
See http://probability.univ.kiev.ua/tims/
\\
The work of F.M. has been carried out in the framework of the activities of the National Group of Mathematical Physics (GNFM, INdAM).



\newpage
\section{Introduction}
 In this paper we outline the relevance of the functions of the Mittag-Leffler type  in renewal processes and in particular in the  thinning theory for long-time behaviour with a generic power law waiting time distribution.
In Section 2 we first recall the definition of a generic renewal process  with the related probability distribution functions. 
In Section 3  we discuss the most celebrated renewal process known as the Poisson process defined by an exponential  probability density function.  Its natural fractional generalization is hence discussed in Section 4 by introducing 
  the so-called renewal process of the Mittag-Leffler type, commonly known as the Fractional Poisson process. 
  Then in Section 5 we consider the thinning theory for  a renewal process with a queue of power law type,  thereby leaning the presentation of Gnedenko and Kovalenko in 1968 pointing out the key role of the Mittag-Leffler function.
Finally,  conclusions  are drown in Section 6.

\section{Essentials of  renewal theory}  

The concept of  {\it renewal process} has been developed as a
stochastic model for describing the class of counting processes
for which the times between successive  events are
independent  identically distributed  ($iid$)
non-negative random variables, obeying a given probability law.
These times are  referred to as waiting times
or inter-arrival times.
For more details see \eg the classical
treatises by
Khintchine \cite{Khintchine_BOOK1960},
Cox \cite{Cox_BOOK1967},
Gnedenko \& Kovalenko \cite{Gnedenko-Kovalenko_BOOK1968},
Feller \cite{Feller_BOOK1971},
and the  most recent books by Ross \cite{Ross_BOOK1997},  by 
Beichelt \cite{Beichelt_BOOK2006}, and by Mitov and Omey \cite{Mitov-Omey_BOOK2014}.
jut to cite the treatises that have mostly attracted our attention.  

For a renewal process having waiting times $T_1,T_2, \dots$, let
$$t_0 = 0\,, \q   t_k= \sum_{j=1}^k T_j\,, \q k \ge 1\,. \eqno(2.1)
    $$
That is $t_1 =T_1$ is the time of the first renewal,
$t_2 = T_1 +T_2$ is the time of the second renewal and so on.
In general $t_k$ denotes the $k$th renewal.

The process is specified if we know the probability law for the waiting
times. In this respect we introduce the
{\it probability density function} ($pdf$)
$\phi(t)$
 and  the   (cumulative) distribution function $\Phi(t)$
so defined:
$$ \phi (t) := \frac{d}{dt} \Phi(t) \,,\q
   \Phi(t) := P \l( T \le t\r) = \int_0^t \phi (t')\, dt'\,.
\eqno(2.2)$$
 When the non-negative random variable  represents
 the lifetime of technical systems, it is common to refer
to $\Phi(t)$ as to the {\it failure probability}
and to
$$ \Psi(t) := P \l(T > t\r) = \int_t^\infty \phi (t')\, dt' = 1-\Phi(t)\,,
\eqno(2.3)$$
as to the {\it survival probability}, because $\Phi(t)$ and $\Psi(t)$ are
the respective probabilities that the system does or does not fail
in $(0, T]$. 
A relevant quantity is the {\it counting function} $N(t)$ that indeed defines the renewal process
as
$$ N(t):= \hbox{max} \l\{k | t_k \le t, \;k = 0, 1, 2, \dots\r\}\,,
\eqno(2.4)$$
that represents the effective number of events before or at instant $t$.
In particular we have $\Psi(t) = P \l(N(t) =0\r)\,.$
Continuing in the general theory we set
$F_1(t) = \Phi(t)$, $f_1(t) = \phi (t)$, and in general
$$ F_k(t) :=  P\l(t_k = T_1+ \dots +T_k \le t \r)\,,
 \; f_k(t) = \frac{d}{dt} F_k(t)\,, \; k\ge 1\,, \eqno(2.5)$$
thus  $F_k(t)$ represents the probability that the sum of the first $k$
waiting times is less or equal  $t$
and $f_k(t)$ its density. Then, for any fixed $k\ge 1$  the normalization
condition  for $F_k(t)$ is fulfilled
because
$$         \lim_{t \to \infty} F_k(t) =
  P\l(t_k = T_1+ \dots +T_k < \infty \r)= 1 \,.
                 \eqno(2.6)$$
In fact,    the sum of $k$ random variables each of which is
    finite with probability 1 is finite with probability 1 itself.
By setting for consistency $F_0(t) \equiv 1$ and $f_0(t) = \delta(t)$,
the Dirac delta  function in the sense of Gel'fand and Shilov \cite{Gelfand-Shilov_1964}
 \footnote{
We find it convenient to recall the {\it formal representation}
 of  this generalized function in $\RR^+\,,$
$$\delta(t) := \frac{t^{-1}}{\Gamma(0)}\,, \q t\ge 0\,.$$},
we also note that for $k \ge 0$ we have
$$ P\l(N(t) =k\r) := P\l(t_k \le t \,,\,t_{k+1} > t\r)
=  \int_0^t f_k(t')\, \Psi (t-t')\, dt'\,.
\eqno(2.7)$$

We now find it convenient to introduce the simplified $\, *\,$ notation
for the  convolution between two causal well-behaved
(generalized) functions $f(t)$ and $g(t)$
$$   \int_0^t   f(t')\, g(t-t')\, dt' = \l(f \,*\, g\r) (t) =
     \l(g \,* \, f\r) (t)  = \int_0^t   f(t-t')\, g(t')\, dt'\,. $$
Being  $f_k(t)$ the $pdf$ of the sum of the  $k$
$iid$ random variables
$T_1,  \dots, T_k$
with $pdf$ $\phi (t)\,, $ we easily recognize that
$f_k(t)$ turns out to be the $k$-fold convolution of $\phi(t)$
with itself,
$$ f_k(t) =  \l(\phi^{*k}\r) (t)\,, \eqno(2.8)$$
so Eq. (2.7)  simply reads:
$$    P\l(N(t) =k\r) = \l(\phi^{*k} \,*\, \Psi\r)(t)\,. \eqno(2.9)$$
Because of the presence of  convolutions
a renewal process is suited for the Laplace transform method.
Throughout this paper we will denote
by   $\widetilde f(s)$
the Laplace transform
of a sufficiently well-behaved (generalized) function
$f(t)$  according to
$$   {\L} \l\{ f(t);s\r\}=   \widetilde f(s)
 = \int_0^{+\infty} \e^{\ds \, -st}\, f(t)\, dt\,,
\q s > s_0\,,
$$
and for $\delta (t)$
consistently we will have
$  \widetilde \delta (s) \equiv 1\,. $
{\it} Note that for our purposes we agree to take $s$ real.
We recognize  that
(2.9) reads in the Laplace domain
$$    \L \{P\l( N(t)=k\r); s\}
= \l[{\widetilde \phi (s)} \r]^k \,   \widetilde \Psi (s)
             \,,\eqno(2.10)$$
where, using (2.3),
$$ \widetilde \Psi (s)  =
\frac{ 1- \widetilde \phi (s)} {s}\,.
\eqno(2.11)$$

\section{The Poisson process as a renewal process}  

The most celebrated   
renewal process is the Poisson process
characterized by a waiting time $pdf$ of exponential type,
$$  \phi (t) = \lambda \, \e^{-\lambda t}\,,\q \lambda >0\,, \q t\ge 0\,.
\eqno(3.1)$$
The process has {\it no memory} being a L\'evy process.
The moments of waiting times of order 1, 2, \dots n  turn out to be
 $$ \langle T\rangle = \frac{1}{\lambda }\,, \q
   \langle T^2 \rangle = \frac{1}{\lambda^2 }\,,
\q \dots \,,  \q \langle T^n \rangle = \frac{1}{\lambda^n }\,, \q \dots \,,
\eqno(3.2)$$
and the {\it survival probability} is
  $$ \Psi(t) := P\l(T>t\r) = \e^{-\lambda t}\,, \q t \ge 0 \,.
\eqno(3.3)$$
We know that the probability  that $k$ events occur in the
interval of length $t $ is
$$ P\l( N(t) = k\r) =   \frac{(\lambda t)^k}{k!} \, \e^{-\lambda t}
\,, \q t \ge 0\,, \q k = 0,1, 2, \dots\;. \eqno(3.4)$$
The probability distribution related to the sum of $k$ $iid$
exponential random variables  is known to be
the so-called {\it Erlang distribution}  (of order $k$).
The corresponding density (the {\it Erlang} $pdf$) is thus
$$ f_k(t) = \lambda \,  \frac{(\lambda t)^{k-1}}{(k-1)!}
    \, \e^{-\lambda t}\,,\q t \ge 0 \,,
\q k =1,2, \dots \,, \eqno(3.5)$$
so that the  Erlang distribution function of order $k$ turns out
to be
$$ F_k(t) = \int _0^t f_k(t')\, dt' =
     1 -  \sum_{n=0}^{k-1}
      \frac{(\lambda t)^n}{n!} \, \e^{-\lambda t}  =
 \sum_{n=k}^{\infty}
      \frac{(\lambda t)^n}{n!} \, \e^{-\lambda t}\,,\q t \ge 0
\,. \eqno(3.6)$$
In the limiting case $k=0$ we recover
$f_0(t) = \delta(t),\; F_0(t) \equiv 1,\; t\ge 0$.

The formulas  (3.4)-(3.6) can easily obtained
by using the technique of the Laplace transform sketched in
the previous section
noting that for the Poisson process we have:
$$ \widetilde\phi(s) =
\frac{\lambda } {\lambda +s}\,,
 \q  \widetilde\Psi(s) =  \frac  {1 } {\lambda +s}\,,   \eqno(3.7)$$
and for the Erlang distribution:
$$  \widetilde f_k(s) =  [\widetilde\phi(s)]^k =
    \frac  {\lambda^k } {(\lambda +s)^k}\,,
\q    \widetilde F_k(s) =  \frac{[\widetilde\phi(s)]^k}{s} =
    \frac  {\lambda^k } {s (\lambda +s)^k}\,.
\eqno(3.8)$$

We also recall that the survival probability
for the Poisson renewal process
obeys the ordinary differential equation  (of relaxation type)
$$     \frac{d}{dt} \Psi(t) = -\lambda \Psi(t)\,, \q t \ge 0\,; \q
\Psi(0^+) =1\,. \eqno(3.9)$$

\section{The renewal  process of the Mittag-Leffler type}   

A "fractional" generalization of the   Poisson renewal process
is simply obtained by generalizing the differential equation  (3.9)
replacing there the first derivative with the integro-differential operator
$\, _tD_*^\beta$ that  is interpreted  as
the fractional derivative
of order $\beta $ in  Caputo's sense.
\vsp
 For a sufficiently well-behaved function $f(t)$ ($t\ge 0$) we
 define the {\it Caputo time fractional derivative}  of
order $\beta  $  with $0<\beta <1$
through
$$
 {\mathcal L} \left\{ _tD_*^\beta \,f(t) ;s\right\} =
      s^\beta \,  \widetilde f(s)
   -s^{\beta  -1}\, f(0^+) 
   \,, \quad f(0^+):= \lim_{t \to 0^+} f(t)\,,
   \eqno(4.1) $$ 
 so that
 $$
    _tD_*^\beta \,f(t) :=
\frac{1}{\Gamma(1-\beta )}\,\int_0^t
 \frac{f^\prime (\tau)}{ (t-\tau )^\beta} \, d\tau\,, \q 0<\beta<1 \,.
 \eqno(4.2) $$
Such operator  has been referred to as
the {\it Caputo} fractional derivative since it
was introduced by Caputo in the late 1960's
for modelling the energy dissipation
 in the rheology of the Earth, see  \cite{Caputo_67,Caputo_69}.
  Soon later this derivative was adopted by Caputo and Mainardi
  in the framework of the linear theory of viscoelasticity, see \cite{CaputoMaina_71}.
\vsp
The reader should observe that the {\it Caputo} fractional derivative
differs from the usual {\it Riemann-Liouville} (R-L) fractional derivative
$$
 \,_tD^\beta  \,f(t) :=
  {\ds {d\over dt}}\,\left[
  {\ds \rec{\Gamma(1-\beta )}\,\int_0^t
    {f(\tau)\,d\tau  \over (t-\tau )^{\beta  }} }\right] \,, \q 0<\beta<1 \,.
\eqno(4.3) $$
\vsp
Following the approach by Mainardi et al.
\cite{Mainardi-Gorenflo-Scalas_VIETNAM04},
we write, taking for simplicity $\lambda =1$,
$$       \, _tD_*^\beta \, \Psi(t) =
- \Psi(t)\,, \q t >0\,,\q 0<\beta \le 1\,; \q \Psi(0^+) =1\,. \eqno(4.4)$$
We also allow
the limiting case $\beta =1$ where all the results of the previous
section (with $\lambda =1$) are expected to be recovered.

We also allow
the limiting case $\beta =1$ where all the results of the previous
sub-section (with $\lambda =1$) are expected to be recovered.
In fact, taking $\lambda=1$ simply means a normalized way of scaling the variable $t$.

For our purpose we need to recall the Mittag-Leffler function
as the natural "fractional" generalization
of the exponential function, that characterizes the Poisson process.
The Mittag-Leffler function of parameter  $\beta\,$
is defined in the complex plane by the power series
$$ E_\beta (z) :=
    \sum_{n=0}^{\infty}\,
   {z^{n}\over\Gamma(\beta\,n+1)}\,, \q \beta >0\,, \q z \in \CC\,.
 \eqno  (4.5)$$
It turns out  to be an entire function of order $\beta $
which reduces for $\beta=1$ to $\exp (z)\,.$
For detailed information on the Mittag-Leffler-type functions
and their Laplace transforms the reader  may  consult \eg
\cite{Erdelyi_HTF,GorMai_CISM97,Podlubny_99} and the most recent monograph by Gorenflo et al.
\cite{Gorenflo-et-al_BOOK2014}.

The solution  of Eq. (4.4) is known to be, see \eg
\cite{CaputoMaina_71,MaiGor_FCAA97}
$$\Psi(t) =  E_\beta (-t^\beta)\,, \q t \ge 0\,, \q 0<\beta \le 1\,,
\eqno (4.6)$$ so
$$ \phi(t) :=    -   \frac{d} {dt} \Psi(t) =
            -   \frac{d}{dt}  E_\beta (-t^\beta) \,, \q t \ge 0
 \,, \q 0<\beta \le 1\,.\eqno (4.7)$$
Then, the corresponding Laplace transforms read
 $$ \widetilde \Psi(s) =
 \frac{s^{\beta-1}}{1+ s^\beta}\,,   \q
 \widetilde \phi(s)= \frac{1} {1+  s^{\beta}}\,,\q
    0<\beta \le 1\,.\eqno(4.8) $$
Hereafter, we find it convenient to summarize
the  most relevant features  of the functions $\Psi(t)$ and $\phi(t)$
when  $0< \beta <1\,.$
We begin to quote their series expansions for $t \to 0^+ $
and  asymptotics for  $t\to \infty $,
$$ \Psi(t)
   = {\ds \sum_{n=0}^{\infty}}\,
  (-1)^n {\ds \frac{t^{\beta n}}{\Gamma(\beta\,n+1)}}
 \,\sim \,  {\ds \frac{\sin \,(\beta \pi)}{\pi}}
  \,{\ds  \frac{\Gamma(\beta)}{ t^\beta}}, \quad  t\to \infty,  \q  0<\beta<1,
     \eqno(4.9) $$
and
$$ \phi(t)
= {\ds \frac{1}{ t^{1-\beta}}}\, {\ds \sum_{n=0}^{\infty}}\,
  (-1)^n {\ds \frac{t^{\beta n}}{\Gamma(\beta\,n+\beta )}}
 \, \sim \,  {\ds \frac{\sin \,(\beta \pi)}{\pi}}
  \,{\ds  \frac{\Gamma(\beta+1)}{ t^{\beta+1}}}\,, \q  t\to \infty, \q 0<\beta <1\,.
     \eqno(4.10) $$
In contrast 
to the Poissonian case  $\beta=1$,
in the case  $0<\beta <1$ for large $t$
the functions $\Psi(t)$ and $\phi(t)$
no longer decay   exponentially  but algebraically.
As a consequence of the power-law asymptotics
the fractional Poisson process for $\beta<1$ turns be no longer Markovian as for $\beta=1$ but of long-memory type.
However, we recognize that for $0<\beta <1$ both  functions
 $\Psi(t)$, $\phi(t)$
keep   the "completely monotonic" character of the Poissonian case.
Complete monotonicity of  the   functions
 $\Psi(t)$ and $\phi(t)$  means
$$ (-1)^n {d^n\over dt^n}\, \Psi  (t) \ge 0\,,  \q
   (-1)^n \frac{d^n}{dt^n}\, \phi  (t) \ge 0\,,
\q n=0,1,2,\dots   \,, \q t \ge 0\,,     \eqno(4.11)$$
or equivalently, their representability as real Laplace transforms
of non-negative   generalized
functions (or measures), see \eg \cite{Feller_BOOK1971}. 

For the generalizations
of Eqs (3.4) and (3.5)-(3.6),  characteristic
of the Poisson and Erlang distributions respectively,
we must point out the   Laplace transform 
$$ \L\{ t^ {\beta \,k}\, E_\beta ^{(k)}
  (-t^\beta ) ;s\} =
        \frac{ k!\, s^{\beta -1}}{(1+s^\beta )^{k+1}}
\,, \q \beta >0 \,, \q k = 0,1, 2, \dots \,,
\eqno(4.12)$$
with $ {\ds E_\beta ^{(k)}(z) := \frac{d^k}{dz^k}  E_\beta(z)}\,, $
that can be deduced  from the book by Podlubny,
see (4.80) in  \cite{Podlubny_99}.
Then, by using the Laplace transforms  (4.8) and
Eqs (4.6), (4.7), (4.12)
in Eqs  (2.8) and (2.9),
 we have  the {\it generalized Poisson distribution},
$$ P\l( N(t) = k\r) =   \frac{ t^{ k\, \beta}}{k!} \,
  E_\beta^{(k)} (- t^\beta)
\,, \q k = 0, 1, 2, \dots \eqno(4.13)$$
and the {\it generalized Erlang} $pdf$   (of order $k \ge 1$),
$$ f_k(t) = \beta \,  \frac{ t^{k\beta-1}}{(k-1)!}
    \, E_\beta^{(k)} (- t^\beta)   
\,. \eqno(4.14)$$
The {\it generalized  Erlang distribution function} turns out
to be
$$ F_k(t) = \int _0^t f_k(t')\, dt' =
     1 -  \sum_{n=0}^{k-1}
      \frac{ t^{n \beta}}{n!} \, E_\beta^{(n)} (- t^\beta)  =
 \sum_{n=k}^{\infty}
      \frac{t^{n\beta}}{n!} \, E_\beta^{(n)} (- t^\beta)
\,. \eqno(4.15)$$


For readers' convenience we conclude this section citing other works dealing with  the so-called fractional Poisson process from different point of view, see \eg Laskin \cite{Las03},
Beghin and Orsingher \cite{BegOrs09}.

\section{The Mittag-Leffler distribution as limit for thinned
 renewal processes} 

Now, we provide, in our notation, an outline of the thinning theory
for renewal processes essentially following the 1968 approach by Gnedenko and Kovalenko
 \cite{Gnedenko-Kovalenko_BOOK1968}.
Examples of  thinning processes are provided in the 
 the 2006 book by 
Beichelt \cite{Beichelt_BOOK2006}, where we read
{\it For instance, a cosmic particle
counter registers only $\alpha$--particles and ignores other types of particles. Or, a
reinsurance company is only interested in claims, the size of which exceeds, say, one
million dollars.} 
\vsp
We must note that other authors, like Sz\'antai
\cite{Szantai_71a,Szantai_71b} speak of {\it rarefaction}
in place of thinning.
\vsp
Let us sketch here the essentials of this theory.
Denoting by $t_n$,\, $n=1,2,3, \dots$
the time instants of events of a renewal process, assuming
$0=t_0<t_1<t_2<t_3<\dots $,
with $i.i.d.$  waiting times $T_1 = t_1\,,\,T_k = t_{k}-t_{k-1}$ for $k\ge 2$,
(generically denoted by T),
{\it thinning} (or {\it rarefaction})
means that for each positive
index  $k$  a decision is made:
the event happening in the instant $t_k$ is deleted with probability $p$
or it is maintained with probability $q=1-p$,  $0<q<1$.
This procedure produces a {\it thinned} or {\it rarefied} renewal process
with fewer events (very few events
if $q$ is near zero, the case of particular interest)
in a moderate span of time.
\vsp
To compensate for this loss  we change the unit of time
so that we still have not very few  but still a moderate number
of events in a moderate span of time.
Such change of the unit of time is
equivalent to  rescaling the waiting time,
multiplying it with a positive factor $\tau $ so that we have
waiting times $\tau T_1,\tau T_2, \tau T_3, \dots$, and
instants $\tau t_1,\tau t_2, \tau t_3,\dots$, in the
rescaled  process.
\vsp 
In other words to bring the distant future into near sight
we   change the unit of time from $1$ to $1/\tau $,
$0 <\tau  \ll 1$.
\vsp
For the random waiting times $T$ this means replacing
$T$ by $\tau T$.
Then, having very many events in a moderate span of time
we compensate this compression by respeeding the whole process,
actually slowing it down so that again we have a moderate
number of events in a moderate span of time.
Our intention is, vaguely speaking, to dispose
on $\tau $ in relation to the rarefaction parameter $q$
in such a way that for $q$ near zero in some sense
the ``average" number of events per unit of time
remains unchanged. In an asymptotic sense
we will make these considerations precise.
\vsp
Denoting by $F(t) = P(T\le t)$ the probability distribution function
of the (original) waiting time $T$,   by $f(t)$ its density
($f(t)$ is a generalized function generating a probability measure)
so that
$F(t) = \int_0^t f(t') \, dt'$,  and analogously by
$F_k(t)$  and $f_k$(t) the distribution 
and density, 
respectively,  of the sum of $k$  waiting times, we have recursively
\begin{equation}
 f_1(t) = f(t) \,,\q
      f_k(t) = \int_0^t  f_{k-1}(t-t') \, dF(t') \,,\;
\hbox{for} \; k \ge 2\,. 
\end{equation}
Observing that after a maintained event the next one of the
original process is kept with probability $q$ but dropped in favor
of the second-next with probability $p\,q$
and, generally, $n-1$ events are dropped  in favor of the
$n$-th-next with probability $p^{n-1}\,q$,
we get for   the waiting time density of the thinned process
the formula
\begin{equation}
g_q(t) = \sum_{n=1}^\infty q\, p^{n-1}\, f_n(t)\,.
\end{equation}
With the modified waiting time $\tau \,T$ we have
$$ P(\tau T\le t) = P(T\le t/\tau ) = F(t/\tau )\,, $$
hence the density $f(t/\tau )/\tau $, and analogously for the
density of the sum of $n$ waiting times
$f_n(t/\tau )/\tau $.
The density of the waiting time of the rescaled (and thinned) process
now turns out as
\begin{equation}
g_{q,\tau} (t) = \sum_{n=1}^\infty q\, p^{n-1}\, f_n(t/\tau)/\tau \,.
\end{equation}
\vsp
In the Laplace domain we have
 $\widetilde f_n(s) = \left(\widetilde f(s)\right)^n\,,$
hence (using $p =1-q$)
\begin{equation}
\widetilde g_q (s)=
\sum_{n=1}^\infty q\, p^{n-1}\,\left(\widetilde f(s)\right)^n
 = \frac{ q\, \widetilde f(s)}{ 1 - (1-q)\, \widetilde f(s)}
 \,,
\end{equation}
from which by Laplace inversion we can, in principle, construct
the waiting time density of the thinned process.
By  re-scaling we get
\begin{equation}
\widetilde g_{q,\tau}(s)=
\sum_{n=1}^\infty q\, p^{n-1}\,\left(\widetilde f(\tau s)\right)^n
 = \frac{ q\, \widetilde f(\tau s)}{ 1 - (1-q)\, \widetilde f(\tau s)}
 \,.
 \end{equation}
\vs
\noindent
     Being interested in stronger and  stronger
thinning ({\it infinite thinning})
let us now  consider 
 a scale of processes with  the parameters $\tau  $ (of {\it rescaling})
and $q$ (of {\it thinning}), with $q$ tending
to zero  {\it under a scaling relation $q = q(\tau) $
  yet to be specified}.
\vsp
We have essentially two cases for the waiting time distribution:
its expectation value (the first moment) is finite or infinite.
In the first case we put
\begin{equation}
\mu = \int_0^\infty t'\,  f(t')\, dt' < \infty \,, \quad \beta=1.
\end{equation}
In the second case we assume a queue of power law type
\begin{equation}
 \Psi(t) :=  \int_t^\infty f(t')\, dt'
\sim \frac{c}{\beta } t^{-\beta }\,, \; t\to \infty, \quad
0<\beta <1\,.
\end{equation}
Then, by the Tauberian theory (see e.g.  \cite{Feller_BOOK1971,Widder_BOOK1946})
the above  conditions mean in the Laplace domain
\begin{equation}
\widetilde f(s) =1- \mu\, s^\beta  + o\left( s^\beta \right)\,,
 \q \hbox{for} \q s \to 0^+\,,
 \end{equation}
with  a positive coefficient $\mu $ and $0<\beta \le 1$.
The case $\beta =1$  obviously corresponds to the situation with finite
first moment, whereas the case $0<\beta <1$ is related
to a power law   queue with
$c= \mu\,\Gamma(\beta +1)\,\sin(\beta \pi)/\pi\,.$
\vsp
Now, passing  to the limit of $q \to 0$ of infinite thinning under the scaling relation
\begin{equation}
 q = \mu \, \tau ^\beta \,, \q 0<\beta \le 1\,, 
\end{equation}
between the positive parameters $q$ and $\tau $,
the Laplace transform of the rescaled density
$\widetilde {g_{q,\tau }}(s)$ in (5.5) of the thinned process tends for fixed $s$ to
\begin{equation}
 \widetilde g(s) = \frac{1}{1+s^\beta}\,,  
 \end{equation}
which corresponds to the Mittag-Leffler density
\begin{equation}
 g(t) = - \frac{d}{dt} E_\beta (-t^\beta )
= \phi^{ML}(t)
\,. 
\end{equation}

We note that in the literature the distribution (5.10) is known as positive Linnik distribution, see e.g. \cite{Devroye_1990}
and used by Pillai \cite{Pillai_1990}to define the Mittag-Leffler distribution. 
Let us remark that Gnedenko and Kovalenko in 1968 obtained
(5.10) as the Laplace transform of the limiting density
but did not identify it as the Laplace transform of a
Mittag-Leffler type function even if this Laplace transform was known since the late 1950's in the Bateman Handbook  \cite{Erdelyi_HTF}.  
Observe once again that in the special case $\beta=1$
we recover  as the limiting process the Poisson process, as formerly shown
in 1956 by R\'enyi \cite{Renyi_1956}.

\section{Conclusions}
We have revised the basic theory of renewal processes including the standard and fractional Poisson processes. and in particular the thinning  theory  with a power law queue,
in order to  out the relevance of the Mittag-Leffler functions.
These processes are essential  in the theory of continuous time random walks,   and, under re-scalng of time and space coordinates, in space-time fractional diffusion processes, as shown e.g. in papers by Gorenflo and Mainardi, see the survey \cite{Gorenflo-Mainardi_KLM2012}.
For the thinning theory we have followed the original approach based on Laplace transform by Gnedenko and Kovalenko, who, however, had ignored the  Mittag-Leffler functions. In practice this note points out a further application of the functions of the Mittag-Leffler type in stochastic processes, not so well known.




\begin{thebibliography}{1}

 \bibitem{BegOrs09}
 L.  Beghin  and E. Orsingher,
  Iterated elastic Brownian motions and  fractional diffusion equations,
 {\it Stochastic Processes and their Applications} {\bf 119} (2009), 1975--2003.



\bibitem{Beichelt_BOOK2006}
F. Beichelt, 
\textit{Stochastic Processes in Science, Engineering and Finance},
Chapman \& Hall/CRC',
Boca Raton, US.
(2006).

\bibitem{Caputo_67}
  {M. Caputo},
  Linear models of dissipation whose $Q$ is almost frequency
  independent,	 Part II.
  {\it Geophys. J. R. Astr. Soc.} {\bf 13}  (1967), 529-539.

\bibitem{Caputo_69}
 {M. Caputo},  {\it Elasticit\`a e Dissipazione}.
 Bologna,  Zanichelli (1969).

 
 \bibitem{CaputoMaina_71}
 {M. Caputo and  F. Mainardi},
  Linear models of  in  anelastic solids.
   {\it Riv. Nuovo Cimento} (Ser. II) {\bf 1} (1971), 161--198.



\bibitem{Cox_BOOK1967}
 D.R. Cox, 
{\it Renewal Theory}, 2-nd Edn, Methuen, London (1967). 

\bibitem{Devroye_1990}
L. Devroye, A note on Linnik's distribution,
{\it  Statistics \& Probability Letters} {\bf 9} (1990), 305--306.
 
\bibitem{Erdelyi_HTF}
  A. Erd\'elyi, W. Magnus, F. Oberhettinger and F.G. Tricomi,
 {\it Higher Transcendental Functions},
  Bateman Project,  
 McGraw-Hill, New York, 1955, Vol 3.
 [Ch. 18: Miscellaneous Functions, pp. 206-227]

\bibitem{Feller_BOOK1971}
{W. Feller},
{\it An Introduction to Probability Theory and its Applications}, Vol. 2,
 2-nd edn. Wiley, New York (1971) .[1-st edn. 1966]

  \bibitem{Gelfand-Shilov_1964} 
  I.M. Gel'fand and G.E.  Shilov, {\it Generalized Functions}, Vol. \textbf{1},
 Academic Press, New York (1964).
 \newline [English translation from the Russian (Nauka, Moscow, 1959)]


\bibitem{Gnedenko-Kovalenko_BOOK1968}
B.V Gnedenko and I.N  Kovalenko, 
\textit{Introduction to Queueing Theory},
Israel Program for Scientific Translations',
Jerusalem (1968).
(Translated from the Russian)

\bibitem{Gorenflo-et-al_BOOK2014}
R. Gorenflo, A. Kilbas, F. Mainardi and  S. Rogosin,
\textit{Mittag-Leffler Functions, Related topics and Applications}
Spinger, Berlin (2014).

\bibitem{GorMai_CISM97}
  R. Gorenflo and F. Mainardi,  Fractional calculus:
  integral and differential equations of fractional order,
  in: A. Carpinteri and F. Mai\-nardi (Editors),
  {\em Fractals and Fractional Calculus in Continuum Mechanics\/},
  Springer Verlag, Wien (1997),  pp. 223--276.

\bibitem{Gorenflo-Mainardi_BAD-HONNEF2008}
R. Gorenflo and  F. Mainardi,
\textit{Anomalous Transport: Foundations and Applications},
Ch 4: Continuous time random walk, Mittag-Leffler waiting time and fractional diffusion: mathematical aspects,
(R. Klages, G. Radons, I.M. Sololov, eds.),
Wiley-VCH,
Weinheim, Germany,
(2008), pp. 93--127.
 [E-print: {\tt http://arxiv.org/abs/0705.0797}]

\bibitem{Gorenflo-Mainardi_KLM2012}  
R.  Gorenflo and F. Mainardi,
Parametric Subordination in Fractional Diffusion Processes, in
 J. Klafter, S.C. Lim  and  R. Metzler  (Editors), {\it Fractional Dynamics},
World Scientific, Singapore (2012), Chapter 10, pp. 229--263
 [E-print: {\tt http://arxiv.org/abs/1210.8414}]

\bibitem{Khintchine_BOOK1960}
 A. Ya. Khintchine, 
{\it Mathematical Methods in the Theory of Queueing}, 
Charles Griffin, London (1960).

\bibitem{Las03} 
N. Laskin,
Fractional Poisson processes,
{\it Comm. Nonlinear Sci. Num. Sim.} {\bf 8} (2003),  201--213.

\bibitem{MaiGor_FCAA97}
F. Mainardi  and R.  Gorenflo,
Time-fractional derivatives in
relaxation processes: a tutorial survey, 
{\it Fract. Calc. Appl. Anal.} \textbf{10} (2007), 269--308.


\bibitem{Mainardi-Gorenflo-Scalas_VIETNAM04}
F. Mainardi, R. Gorenflo and  E. Scalas,
    \textit{A fractional generalization of the Poisson processes}
    {Vietnam Journal of Mathematics} {\bf 32} SI (2004), 53--64.
	

\bibitem{Mitov-Omey_BOOK2014}
K.V. Mitov and E. Omey,
{\it Renewal Processes},
Springer, Heidelberg (2014).

\bibitem{Pillai_1990}
R.N. Pillai, On Mittag-Leffler functions and related distributions,
{\it Ann. Inst. Statist. Math.} {\bf 42} No 1 (1990), 157--161.

\bibitem{Podlubny_99} 
I. Podlubny,  {\it Fractional Differential Equations},
Academic Press, New York (1999).

\bibitem{Renyi_1956}
A. Renyi,
 \textit{A characteristic of the Poisson stream},
 { Proc. Math. Inst. Hungarica Acad. Sci.}  {\bf 1} (1956), no. 4, 563--570. 
 [In Hungarian]
 
 \bibitem{Ross_BOOK1997} 
  S.M. Ross,  {\it Introduction to Probability Models},
6-th Edn, Academic Press, New York (1997).
 
\bibitem{Szantai_71a}
T. Sz\`antai,
\textit{Limiting distribution for the sums
of random number of random variables concerning the
rarefaction of recurrent events},
{Studia Scientiarum Mathematicarum Hungarica} {\bf 6} (1971), 443--452.

\bibitem{Szantai_71b}
T. Sz\`antai,
\textit{On an invariance problem related to
different rarefactions of recurrent event},
{Studia Scientiarum Mathematicarum Hungarica} {\bf 6} (1971), 453--456.

\bibitem{Widder_BOOK1946}
 D.V. Widder,
 {\it The Laplace Transform},
 Princeton University  Press, Princeton (1946).

\end{thebibliography}
\end{document}